\documentstyle{article}\begin{document}\centerline{\bf A sum rule for the critical zeros of $\zeta(s)$}\vskip .3in

\centerline{ M.L. Glasser}

\centerline{Department of physics, Clarkson University}
 
\centerline{Potsdam, NY 13699-5820}
\centerline{laryg@clarkson.edu}
\vskip .5in

\centerline{\bf Abstract}\vskip .1in
\begin{quote}
The two-parameter  series over the critical zeros of the Riemann Zeta function 

$$Re\sum_{\rho}\frac{ x^{(\rho-a)/4a}}{\sqrt{\rho-a}\sinh[\frac{\pi}{2}\sqrt{\frac{\rho-a}{a}}]\zeta'(\rho)}$$
is evaluated in terms of $\zeta(s)$  on the real axis.
\end{quote}
 
\newpage

\section{Introduction}

Recently, Guillera[1] has derived the important sum rule for the nontrivial zeros, denoted $\rho$, of the Riemann zeta function,
$$\sum_{\rho}\frac{x^{\rho-1/2}}{\sin\pi(\rho-1/2)}=\sqrt{x}-\frac{1}{\pi}\frac{\zeta'(1/2)}{\zeta(1/2)}+h(x)+\frac{1-x^2}{\pi}\sum_{n=1}^{\infty}\frac{\sqrt{n}\Lambda(n)}{(n+x)(1+nx)}\eqno(0.1)$$
$$h(x)=\frac{1}{\sqrt{x}(x^2-1)}-\frac{1}{2x-2}+\frac{\ln(8\pi)+\gamma}{\pi(x+1)}-\frac{2}{\pi}\frac{\sqrt{x}\cot^{-1}\sqrt{x}}{x+1},$$
where $\Lambda$ denotes Mangolt's function[2].  Guillera  then 
demonstrated it's utility, both assuming the Riemann Hypothesis (RH) and in general, for elucidating  additional  properties of the zeros and consequences of the RH. This suggests that if (0.1) could be generalized to include an additional parameter, that perhaps further information could be obtained. Unfortunately, an attempt to achieve this failed, but a somewhat more complicated relative of Guillera's relation was found involving a second free parameter and presenting the, possibly advantageous, feature that the Mangoldt and Chebyshev functions do not appear. This note is concerned only with a brief derivation, but further series and their possible consequences will be presented elsewhere.

\section{Calculation}

We start with the integral[3]
$$\int_{-i\infty}^{i\infty}\frac{ds}{2\pi i}\frac{x^{s(1-s)}}{\cos(\pi s)\zeta[4a s(1-s)]}=\frac{x^{1/4}}{2\pi\zeta(a)},\eqno(1.1)$$
where $a>0$ and $0<x<1$.
The singularities of the integrand in the right half $s-$plane are the simple poles
$$s_n=\frac{1}{2}[1+\sqrt{\frac{2n+a}{a}}],\qquad n=1,2,3\cdots.\eqno(1.2a)$$
$$s_k=k+\frac{1}{2},\qquad k=0,1,2,\cdots.\eqno(1.2b)$$
$$s_{\rho}=\frac{1}{2}[1+\sqrt{\rho/a-1}]\mbox{ and }s_{\rho^*} \eqno(1.2c)$$
where $\{\rho\}$ is the set of critical zeros in the upper half plane. 
Since the integrand of (1.1) decays exponentially into the right half $s-$plane, the contour can be closed by the right half of an infinite circle and we obtain by summing the corresponding residues
$$Re\sum_{\rho}\frac{ x^{(\rho-a)/4a}}{\sqrt{a-\rho}\sin[\frac{\pi}{2\sqrt{a}}\sqrt{a-\rho}]\zeta'(\rho)}$$
$$=\frac{\sqrt{a}}{\pi\zeta(a)}+\sum_{n=1}^{\infty}\frac{(-1)^{n+1}(2\pi)^{2n}x^{-(2n+a)/4a}}{\sqrt{2n+a}\sin[\frac{\pi}{2}\sqrt{\frac{2n+a}{a}}]\zeta(2n+1)(2n)!}-\frac{2\sqrt{a}}{\pi}\sum_{k=1}^{\infty}\frac{(-1)^k x^{-k^2}}{\zeta[a(1-4k^2)]},\eqno(1.3)$$
where we have used
$$\zeta'(-2n)=\frac{(-1)^n\zeta(2n+1)(2n)!}{2(2\pi)^{2n}}.\eqno(1.4)$$
From (1.3) one sees that the critical zeros are related to the values of $\zeta$ at certain points on the real axis.  It appears that  there are many analogues of (1.1) that lead to similar sum rules, but with more free parameters. It is possible that by manipulating these, interesting information might be obtained concerning  the zeros. This will be examined in a subsequent report. From (1.3) with $a=1/2$ we find that the RH is equivalent to
$$Re\sum_{\tau>0}\frac{\tau^{-1/2} e^{\frac{1}{2}i(\tau\ln x+\frac{\pi}{2})}}{\sin\left[\frac{\pi\sqrt{\tau}}{1+i}\right]\zeta'(\frac{1}{2}+i\tau)}=\frac{1}{\pi\sqrt{2}\zeta(1/2)}
-\frac{\sqrt{2}}{\pi}x^{1/4}\sum_{k=1}^{\infty}\frac{(-1)^{k+1}}{x^{k^2}\zeta(\frac{1}{2}-2k^2)}$$
$$+\sqrt{2}x^{-1/4}\sum_{n=1}^{\infty}\frac{(-1)^{n+1}(2\pi)^{2n}}{x^n\sqrt{4n+1}\sin\left(\frac{\pi}{2}\sqrt{4n+1}\right)\zeta(2n+1)(2n)!}\eqno(1.5)$$
where $\{\tau\}$ is the set of imaginary parts of the critical zeros.

\vskip .2in
\noindent
{\bf Acknowledgements} The author thanks Dr. Michael Milgram for helpful correspondence.

\noindent
{\bf References}
\vskip .1in
\noindent
[1] J. Guillera,  arXiv: 1307.5723v1

\noindent
[2] E.C. Tichmarsh,{\it The Riemann Zeta Function}, [Cambridge Univ. Press, London 1930]

\noindent
[3] M.L. Glasser, arXiv:1308.6361v2

\end{document}